\documentclass[a4paper,11pt,onecolumn]{article}

\usepackage{bbm}
\usepackage{amsfonts}

\usepackage{stmaryrd}

\usepackage{graphicx}

\usepackage{geometry}
\usepackage{fancyhdr}
\usepackage{fancyvrb}
\usepackage{fancybox}
\usepackage{amsmath,amsfonts,amssymb,amscd,graphicx}
\usepackage{subfigure}
\usepackage{indentfirst}
\usepackage{bm}
\usepackage{multicol}
\usepackage{abstract}
\usepackage{anysize}
\usepackage{hyperref}
\usepackage{listings}
\usepackage{color}
\usepackage{units}
\usepackage{tabularx}
\usepackage{mathrsfs}
\usepackage[all]{xy}

\usepackage{extarrows}

       \font\tenmsb=msbm10
       \font\sevenmsb=msbm7
       \font\fivemsb=msbm5
       \catcode`\@=11
       \ifx\amstexloaded@\relax\catcode`\@=\active
       \endinput\else\let\amstexloaded@\relax\fi
       \def\spaces@{\space\space\space\space\space}
       \def\spaces@@{\spaces@\spaces@\spaces@\spaces@\spaces@}
       \def\space@.{\futurelet\space@\relax}
       \space@. %
       \def\Err@#1{\errhelp\defaulthelp@\errmessage{AmS-TeX error: #1}}
       \def\relaxnext@{\let\next\relax}
       \def\accentfam@{7}
       \def\noaccents@{\def\accentfam@{0}}

       \def\Cal{\relaxnext@\ifmmode\let\next\Cal@\else
       \def\next{\Err@{Use \string\Cal\space only in math mode}}\fi\next}
       \def\Cal@#1{{\Cal@@{#1}}}
       \def\Cal@@#1{\noaccents@\fam\tw@#1}

       \def\Bbb{\relaxnext@\ifmmode\let\next\Bbb@\else
       \def\next{\Err@{Use \string\Bbb\space only in math mode}}\fi\next}
       \def\Bbb@#1{{\Bbb@@{#1}}}
       \def\Bbb@@#1{\noaccents@\fam\msbfam#1}

       \newfam\msbfam
       \textfont\msbfam=\tenmsb
       \scriptfont\msbfam=\sevenmsb
       \scriptscriptfont\msbfam=\fivemsb

\def\Z{\mathbb{Z}}
\def\R{\mathbb{R}}
\def\C{\mathbb{C}}

\def\g{\mathfrak{g}}

\def\i{\mathrm{i}}

\newtheorem{theorem}{Theorem}
\newtheorem{lemma}{Lemma}[section]
\newtheorem{proposition}{Proposition}[section]
\newtheorem{corollary}[proposition]{Corollary}

\newtheorem{assumption}{Assumption}

\newtheorem{remark}{Remark}

\newtheorem{example}{Example}
\newtheorem{definition}{Definition}[section]

\@addtoreset{equation}{section}

\newcommand{\qed}{\nolinebreak\hfill$\Box$
\par\medbreak}

\newcommand{\la}{\langle}
\newcommand{\ra}{\rangle}

\newcommand{\bt}{\begin{theorem}}
\newcommand{\et}{\end{theorem}}

\newcommand{\bq}{\begin{equation}}
\newcommand{\eq}{\end{equation}}

\newcommand{\bp}{\begin{proposition}}
\newcommand{\ep}{\end{proposition}}

\newcommand{\bc}{\begin{corollary}}
\newcommand{\ec}{\end{corollary}}

\newcommand{\bdf}{\begin{definition}}
\newcommand{\edf}{\end{definition}}

\newcommand{\bl}{\begin{lemma}}
\newcommand{\el}{\end{lemma}}

\newcommand{\bay}{\begin{array}}
\newcommand{\ea}{\end{array}}

\newcommand{\ba}{\begin{eqnarray}}
\newcommand{\na}{\end{eqnarray}}

\newcommand{\bas}{\begin{assumption}}
\newcommand{\eas}{\end{assumption}}

\newcommand{\bx}{\begin{example}}
\newcommand{\ex}{\end{example}}

\newcommand{\br}{\begin{remark}}
\newcommand{\er}{\end{remark}}

  \newcommand{\cE}{{\mathcal E}}

  \newcommand{\cP}{{\mathcal P}}

  \newcommand{\cO}{{\mathcal O }}

     \newcommand{\cS}{{\mathcal  S}}

  \newcommand{\D}{{\mathbb D}}

\begin{document}
\setlength{\columnsep}{11pt}
\title{Sum of Hamiltonian manifolds}
\author{\small  Bohui Chen, \   Hai-Long Her\  \   and    \     Bai-Ling Wang}

\date{}
\maketitle

\begin{quote}
\small {{\bf Abstract}.\ 
 For any compact connected Lie group $G$,   we study the Hamiltonian  sum of  two compact  Hamiltonian group $G$-manifolds $(X^+,\omega^+,\mu^+)$ and $(X^-,\omega^-,\mu^-)$     with a common codimension 2 Hamiltonian submanifold  $Z$ of the  opposite   equivariant Euler classes of the normal bundles.    We  establish   that  the symplectic reduction of the Hamiltonian sum agrees with the symplectic sum of 
the  reduced symplectic manifolds. We also compare the equivariant first Chern class of  the Hamiltonian sum  with   the equivariant first Chern classes of  
$X^\pm$.  }
\end{quote}




\tableofcontents

\section{Introduction}\label{SEC-1}

Symplectic cut \cite{L1} and symplectic  sum \cite{G, MW} are two of important operations on symplectic manifolds. 
These are  effective methods of constructing new symplectic manifolds with specific properties, which can apply to
study the problem how some important symplectic invariants, such as Gromov-Witten (GW) invariants, 
change under such symplectic operations (\cite{T,LR,IP2,H,FZ}). 
Li and Ruan \cite{LR} studied the degeneration formula of GW invariants under symplectic cut.
Ionel and Parker \cite{IP2} gave the symplectic sum formula of GW invariants under symplectic sum of two symplectic manifolds 
with a common codimension 2 symplectic submanifold. These  methods were complemented and strengthened by  Tehrani and Zinger \cite{FZ}.

In recent years, motivated by the Yang-Mills-Higgs theory in physics as well as some observations from the proof of Atiyah-Floer conjecture,  
some interesting invariants, called {\it Hamiltonian} or {\it gauged} GW invariants, were discovered for a symplectic manifold $(X,\omega)$ with a Hamiltonian action of a compact Lie group $G$ (\cite{CGS,M,CGMS,MT,CWW}). 
A  symplectic manifold $(X,\omega)$ with a Hamiltonian $G$-action is  called a Hamiltonian $G$-manifold.
The crucial notion for Hamiltonian manifolds is the moment map which  is a $G$-equivariant smooth map 
$\mu:X\to\g^*$ generating the Hamiltonian $G$-action, where $\g^*$ is the dual of the Lie algebra $\g:= \mathrm{Lie}(G)$.
We denote a Hamiltonian $G$-manifold by $(X,\omega, \mu)$. Hamiltonian GW invariants can be constructed by studying the equations of 
symplectic vortices derived from a principal $G$-bundle over a Riemann surface and $(X,\omega, \mu)$.
In this paper, we study the  sum operations on Hamiltonian  $G$-manifolds. 
As a potential application, such sum operation would be useful for 
establishing the  sum formula or degeneration formula  of the  Hamiltonian GW invariants.
 \bigskip

Let $(X^+,\omega^+,\mu^+)$ and $(X^-,\omega^-,\mu^-)$ be two Hamiltonian $G$-manifolds of the same dimension,
$Z=Z^+ \cong Z^-$ be their common codimension 2 Hamiltonian $G$-submanifold {which has the same $G$-action  induced from  both sides}, 
with restricted moment maps, denoted by $\mu_Z=\mu^+|Z=\mu^-|Z$. 
Suppose that the normal bundles $N_{X^+}Z$ and $N_{X^-}Z$  have opposite equivariant Euler classes
$e_G(N_{X^+}Z)=-e_G(N_{X^-}Z)$.
Then we can fix a $G$-equivariant isomorphism of the two $G$-trivial complex line bundles
\bq\label{iso-1}
\Phi:N_{X^+} Z \otimes_\C N_{X^-} Z \xlongrightarrow[\ \ G\ \  ]{\cong} Z\times\C.
\eq
The  symplectic sum  of $(X^+,\omega^+ )$ and $(X^-,\omega^- )$  gives rises  to  
 a $(2n+2)$-dimensional  symplectic manifold $(\cS , \Omega)$, unique up to the (non-equivariant) homotopy class of $\Phi$, a smooth map $\pi : \cS \to D$ where $D$ is  a  sufficiently small neighbourhood of the origin in   $ \C$ such that 
 \begin{itemize}
\item  $\pi$ is surjective and  $\cS_0 =\pi^{-1}(0) = X^+\cup_ZX^-$, 
\item $\pi$ is submersion away from $Z\subset \cS_0$,
\item the restriction of $\Omega$ to $\cS_\lambda = \pi^{-1}(\lambda)$ is nondegenerate for every $\lambda\in D\setminus \{0\}$,
\item $\Omega |_{X^\pm} = \omega^\pm$ for $X^\pm \subset  \cS_0$.  
\end{itemize}

Our main result is the following Hamiltonian  sum theorem which is in principle a combination of 
the symplectic sum construction  in \cite{IP2,FZ} and a gluing construction of moment maps.

\begin{theorem} \label{thm-1}
Let  $G$ be a connected compact Lie group,
$(X^+,\omega^+,  \mu^+)$ and     $(X^-,\omega^-,  \mu^-)$ be two $2n$-dimensional Hamiltonian $G$-manifolds with a 
common codimension $2$ Hamiltonian  submanifold $Z=Z^+ \cong Z^-$  such that  $\mu_Z= \mu^+|_Z=\mu^-|_Z$ and  their respective normal bundles have opposite $G$-equivariant Euler classes. 
Then for each choice of  homotopy class of  equivariant  isomorphisms (\ref{iso-1}), there exist a  natural Hamiltonian structure on   $(\cS , \Omega)$ 
with a moment map $\mu:\cS\to\g^*$,
such that  for every $\lambda\in D\setminus\{0\}$, the fiber $(\cS_\lambda=\pi^{-1}(\lambda),\Omega|_{\cS_\lambda}, \mu_\lambda)$  is also a Hamiltonian $G$-manifold 
with the  moment map  $\mu_\lambda$ given by the restriction of  $\mu$ to  $\cS_\lambda$ and 
$\mu|_{X^\pm} = \mu^\pm$ for $X^\pm \subset \cS_0$.  Moreover, if  $0$ is a regular value of $\mu_Z$ and $G$-actions on 
$(\mu^\pm)^{-1}(0)$ are free, then for every $\lambda\in D\setminus\{0\}$ 
\[
\mu^{-1}_\lambda (0)/G
\]
is the symplectic sum of symplectic reductions  $(\mu^+)^{-1}(0)/G$ and  $(\mu^+)^{-1}(0)/G$ 
along their common codimension 2 symplectic submanifold $\mu_Z^{-1}(0)/G$.  
\end{theorem}
For every $\lambda\in D\setminus\{0\}$, the fiber $(\cS_\lambda, \Omega|_{\cS_\lambda}, \mu_\lambda)$ is  called a 
Hamiltonian sum of  $(X^+,\omega^+,  \mu^+)$ and $(X^-,\omega^-,  \mu^-)$, or simply denoted by $\cS_\lambda=X^+\#^G_{Z,\lambda} X^-$. 
All $\cS_\lambda$ are deformations, in the category of Hamiltonian manifolds, of the singular fibre $\cS_0$. 
For $\lambda\neq 0$, these are  isotopic to one another, as Hamiltonian $G$-manifolds. 

The  symplectic sum of symplectic orbifolds along a symplectic normalisable orbifolds was  constructed  in \cite{VM}.  If  $G$-actions  on $(\mu^\pm)^{-1}(0)$ are locally free, 
symplectic reductions  $(\mu^\pm)^{-1}(0)/G$   are symplectic orbifolds, and  $\mu_Z^{-1}(0)/G$ is a 
symplectic normalisable suborbifold as in \cite{VM}.  Then  the principle that symplectic sum commutes with reduction  holds
in this Hamiltonian case.  
\bigskip

The paper is organized as follows.
In section \ref{SEC-MMP}, we review some  definitions and properties of   Hamiltonian manifolds and moment maps. We can generalize the moment map to  a $G$-action on a smooth manifold with a $G$-invariant closed 2-form.
In section \ref{Sec-rev}, we review   the construction of symplectic  sum   based on the paper \cite{FZ}.     In section \ref{SEC-Sum-Ham}, we show that the symplectic sum of two Hamiltonian $G$-manifolds has a natural Hamiltonian structure.
In section \ref{equi-1-Chern}, we  compare  the equivariant first Chern class of $(\cS_\lambda,  \Omega_\lambda, \mu_\lambda)$ with  
  the equivariant first Chern classes of $(X^+,\omega^+,  \mu^+)$ and     $(X^-,\omega^-,  \mu^-)$.



\section {Hamiltonian manifolds, moment maps and generalization}\label{SEC-MMP}

In this section, we first review some basic notions and an example related to moment maps of Hamiltonian group $G$-action on symplectic manifolds (\cite{MS}), 
then study its generalization to moment maps corresponding to any $G$-invariant closed 2-forms on general manifolds.

\subsection{Moment maps for Hamiltonian $G$-manifolds}

We first consider a symplectic manifold $(X, \omega)$.
Let $G$ be a compact Lie group  which acts covariantly on  $(X, \omega)$ by symplectomorphisms. 
This means that there is a smooth group homomorphism from Lie group $G$ to the group of 
symplectomorphisms $\mathrm{Symp}(X, \omega)$: $g\mapsto\phi_g$.
Denote the Lie algebra of $G$ by $\g:= \mathrm{Lie}(G)$ as right  invariant vector fields on $G$. The infinitesimal action determines a Lie algebra 
homomorphism from $\g$ to the Lie algebra of symplectic vector fileds
$\mathcal{X}(X,\omega):\xi\mapsto V_\xi$ defined by  
\bq\label{Vxi}
V_\xi(x):=\frac{d}{dt}\Bigg|_{t=0}\phi_{\exp(t\xi)}x,\ \ \ \  \forall\ x\in X.
\eq
By  calculation one can verify that
\bq\label{Equiv-Vxi}
V_{\mathrm{Ad}(g^{-1})\xi} = \phi_g^*V_\xi, \ \ \ \  \ \    V_{[\xi,\eta]} = [V_\xi , V_\eta].
\eq
So $V_\xi$ is $G$-equivariant.
A vector field $V_\xi$ is Hamiltonian if there is a corresponding Hamiltonian function $H_\xi$ such that  $\iota(V_\xi)\omega=dH_\xi$. 
The action of $G$ on $X$ is called {\bf Hamiltonian} if the vector field $V_\xi$ is Hamiltonian  for every $\xi\in\g$ and the map
\begin{eqnarray}
\g & \rightarrow& C^\infty(X,\R),\nonumber\\
\xi &\mapsto    &     H_\xi\nonumber
\end{eqnarray}
can be chosen to be $G$-equivariant with respect to the adjoint action of $G$ on its Lie algebra $\g$, that is for $\forall$ $x\in X$,
$$H_{\mathrm{Ad}(g^{-1})\xi} (x)=H_\xi(g\cdot x).$$
A Hamiltonian $G$-manifold is  a  symplectic manifold $(X,\omega)$ with  a Hamiltonian group 
$G$-action. By definition, the $G$-action is generated by the Hamiltonian vector fields $V_\xi$ associated to the Hamiltonian functions
$H_\xi:X\to\R$ such that $\xi\mapsto H_\xi$ is $G$-equivariant.

Then  the (symplectic) moment map for the $G$-action on $(X, \omega)$ is a $G$-equivariant smooth map 
$$\mu:X\to\g^*\cong\g$$ 
(isomorphism via an invariant inner product on the Lie algebra $\g$), such that the Hamiltonian vector fields associated to the Hamiltonian functions defined by 
\bq\label{H-mmp}
H_\xi(x)=\la \mu(x),\xi\ra
\eq
generate the action, where $\la \cdot,\cdot\ra$ denotes the pairing between $\g^*$ and $\g$, or the invariant inner product on $\g$. 
In other words, $\mu$ is a moment map for the $G$-action  if and only if $\mu$ satisfies the following
\begin{itemize}
\item[\rm{(a)}] $G$-equivariance condition 
\bq\label{mu-Adg}
\mu(g\cdot x)=\mathrm{Ad}(g^{-1})^*\mu(x);
\eq
\item[\rm{(b)}] Hamiltonian condition
\bq\label{MMP}
\la d\mu(x)\tau,\xi\ra=\omega(V_\xi(x),\tau)
\eq
\end{itemize}
for all $g\in G$, $\xi\in \g$, $x\in X$, $\tau\in T_xX$. 
We simply say that the Hamiltonian $G$-action is generated by the moment map $\mu$,
and denote a  Hamiltonian G-manifold  by the tuple $(X,\omega, G, \mu)$ or simply $(X,\omega, \mu)$ if the group $G$ is prescribed.
Here is an example of Hamiltonian manifolds.

\bx[Induced Hamiltonian action on cotangent bundle]\label{Exp-mmp3}
In classical mechanics, any cotangent bundle is a symplectic manifold corresponding to the phase space of a Hamiltonian system.
Let $G$ be a Lie group, $M$  a differential manifold with a $G$-action. 
Then there exist associated natural vector bundle $G$-actions on the tangent bundle $TM$ and the cotangent bundle $T^*M$.
The action on $TM$ is given by $g\cdot v = dg_x(v)$, for $v\in T_xM$.
The action on the cotangent bundle $T^*M$ is given by 
$$\la g\cdot\eta,  v\ra_d=\la\eta, g^{-1}\cdot v\ra_d$$ 
for $\forall$ $v\in T_xM$ and $\eta\in T_x^*M$, where $\la\cdot, \cdot \ra_d$ is the dual between $T^*M$ and $TM$.
It is well-known that there exist a canonical 1-form $\theta$ on $T^*M$ defined as
$$\theta_{(x,\alpha)}(v,\eta)=\la\alpha, v \ra,$$
where $(x,\alpha)\in T^*_xM$, $(v,\eta)\in T_{(x,\alpha)}(T^*M)$.
The canonical symplectic form on cotangent bundle $T^*M$ is $$\omega_{can}= - d\theta.$$ 
For such {\bf induced $G$-action} on $T^*M$, the moment map is given by the canonical 1-form
$$\la \mu_{\theta}(\cdot), \xi \ra=-\theta(V_\xi(\cdot) ),$$
for $\forall$  $\xi\in\g$. Condition (\ref{MMP}) obviously holds. Condition (\ref{mu-Adg}) holds because $V_\xi$ is $G$-equivariant.
So $(T^*M, \omega_{can}, G, \mu_{\theta})$ is  a Hamiltonian $G$-manifold. \qed
\ex

\subsection{Generalization of moment maps}\label{OA-mmp}

In fact, for any manifold $M$ with a fixed closed 2-form $\omega$ (might not be non-degenerate), 
one can generalize the notion of moment maps for any group action on $M$ preserving the closed 2-form $\omega$.
Now the pair $(M, \omega)$ might not be a symplectic manifold.
Let $G$ be a compact Lie group  which acts  on the pair $(M, \omega)$ preserving the closed 2-form $\omega$. 
This means that there is a smooth group homomorphism from Lie group $G$ to the group of 
$\omega$-preserving diffeomorphisms $\mathrm{Diff}(M, \omega)$.

For an element $\xi\in\g$, the infinitesimal action determines a vector field $V_\xi$  defined as in (\ref{Vxi}).
Since the closed 2-form $\omega$ is $G$-invariant, by Cartan formula
$$0=\mathcal{L}_{V_\xi}\omega=d(\iota_{V_\xi} \omega)+\iota_{V_\xi}d\omega=d(\iota_{V_\xi} \omega),$$ 
we see that  $\iota_{V_\xi} \omega$ is a closed 1-form.
A smooth map 
\bq\label{muo}
\mu_\omega: M\to\g^*
\eq
satisfying (\ref{mu-Adg}) and (\ref{MMP}) is called the {\bf $\omega$-moment map} for the $G$-action on $(M, \omega)$. 
Equivalently, $\mu_\omega$ is $G$-equivariant and the following equality holds for all $\xi\in \g$
\bq\label{ommp}
d\la\mu_\omega,\xi\ra=\iota_{V_\xi} \omega.
\eq
In particular, if the closed 2-form $\omega$ is non-degenerate, an $\omega$-moment map is just the (symplectic) moment map for Hamiltonian $G$-action on a symplectic manifold.

If the $G$-invariant 2-form $\omega$ is the exterior differential of a $G$-invariant 1-form $\alpha$, that is $\omega=d\alpha$, then by Cartan formula
$$0=\mathcal{L}_{V_\xi}\alpha=d(\iota_{V_\xi} \alpha)+\iota_{V_\xi}d\alpha=d(\iota_{V_\xi} \alpha)+\iota_{V_\xi}\omega.$$
From (\ref{ommp}), this means 
$$d\la\mu_\omega,\xi\ra= - d(\iota_{V_\xi} \alpha).$$
Then a $G$-equivariant smooth map 
\bq\label{mu-alpha}
\mu_\alpha: M\to\g^*
\eq
is called the {\bf $\alpha$-moment map} for the $G$-action on $(M, \omega=d\alpha)$ if the following equality holds for all $\xi\in \g$
\bq\label{al-mmp}
\la\mu_\alpha,\xi\ra=-\iota_{V_\xi} \alpha.
\eq
We remark that  the {\bf $\alpha$-moment map} was first defined by Lerman \cite{L2} and further studied by Chiang and Karshon in \cite{CK} for contact manifolds.   Note that in this case the $d\alpha$-moment map is just the $\alpha$-moment map $\mu_{d\alpha}=\mu_\alpha$. 
 
Let $f:M\to\R$ be a smooth $G$-invariant function. Then it is easy to get the following relation between the $\alpha$-moment map and $f\alpha$-moment map 
\bq\label{rel-alpha}
\mu_{f\alpha}=f\mu_\alpha.
\eq

In  section \ref{SEC-Sum-Ham}, we will consider a  Hamiltonian $G$-manifold $(X,\omega,\mu)$ such that
the symplectic form can be expressed as $\omega=\omega_0+d\alpha$, where $\omega_0$ is a $G$-invariant closed 2-form  and $\alpha$ is a $G$-invariant 1-form on $X$.
The following lemma is obvious. 
\bl\label{Cor-mmp}
Given a symplectic manifold $(X,\omega)$ with a $G$-action. Suppose that the symplectic form $\omega=\omega_0+d\alpha$ such that 
$\omega_0$ is a $G$-invariant closed 2-form  and $\alpha$ is a $G$-invariant  1-form. 
Let $\mu_{\omega_0}$ and $\mu_\alpha$ be the $\omega_0$-moment map and $\alpha$-moment map on $X$. 
Then the sum
\bq\label{mu-sum}
\mu:=\mu_\omega=\mu_{\omega_0}+\mu_\alpha
\eq
is the moment map of the $G$-action on $X$. So $(X,\omega,\mu)$ is a Hamiltonian $G$-manifold. 
\el

Given two manifolds $M$ and $N$ with $G$-action, we consider a smooth $G$-equivariant  map $\varphi:M\to N$.
Let $\omega$ and $\alpha$ be a $G$-invariant closed 2-form and a $G$-invariant 1-form on $N$, respectively. 
Then $\varphi^*\omega$ and $\varphi^*\alpha$ are $G$-invariant closed 2-form and $G$-invariant 1-form on $M$. 
Furthermore, if $\mu_\omega$ and $\mu_\alpha$ are $\omega$-moment map and $\alpha$-moment map on $N$, 
then the $\varphi^*\omega$-moment map and $\varphi^*\alpha$-moment map on $M$ are
\bq\label{phi-mmp}
\mu_{\varphi^*\omega}=\mu_{\omega}\circ\varphi,\ \ \   \mathrm{and}   \ \ \  \mu_{\varphi^*\alpha}=\mu_{\alpha}\circ\varphi,
\eq
respectively.

\section{Review of symplectic sum construction}\label{Sec-rev}

In this section, we  review  the  construction of symplectic sum  based on \cite{FZ}(see also \cite{G, IP2}).
 Let $(X^+,\omega^+)$ and $(X^-,\omega^-)$ be two compact  symplectic manifolds of dimension $2n$,
$Z=Z^+ \cong Z^-$ be their common codimension 2 symplectic submanifold. 
Take an  almost complex structure $J^\pm$ on $X^\pm$ compatible with $\omega^\pm$ such that $J^+|_Z=J^-|_Z$.
Since the normal bundle $N_\pm Z:=N_{X^\pm}Z$  inherits a symplectic structure from $\omega^\pm$ and thus a complex structure up to homotopy, 
both $N_{+}Z$ and $N_{-}Z$ are oriented.  We assume that they have opposite  Euler classes:
\bq\label{5EL}
e(N_{+}Z)+e(N_{-}Z)=0.
\eq
There exists an somorphism  of  complex line  bundles 
\bq\label{5iso-LB}
\Phi:N_{+} Z \otimes_\C N_{-} Z \longrightarrow Z\times\C.
\eq
These data determine a family of symplectic sums. 
The following theorem is taken from \cite[Proposition 3.1]{FZ}.

\bt[Symplectic connect sum]\label{SCS} 
Given two compact $2n$-dimensional symplectic manifolds $(X^+,\omega^+)$ and $(X^-,\omega^-)$ with a 
common codimension 2 symplectic submanifold $Z$ satisfying (\ref{5EL}). 
Then for each choice of homotopy class of  isomorphisms (\ref{5iso-LB}), 
there exist a $(2n+2)$-dimensional symplectic manifold $(\cS , \Omega)$, 
a smooth map  $\pi : \cS \to D$ over a sufficiently small  neighbourhood  $D$ of the origin in  $\C$, and an $\Omega$-compatible almost complex structure $J_{\cS}$ on $\cS$ such that
\begin{itemize}
\item[\rm{(1)}] $\pi$ is surjective and  is a submersion outside of $Z\subset \cS_0$, and $\cS_0=X^+\cup_Z X^-$; 
\item[\rm{(2)}] the restriction $\omega_\lambda$ of  $\Omega$ to $\cS_\lambda=\pi^{-1}(\lambda)$ is nondegenerate for $\forall$ $\lambda\in D^*=D\setminus\{0\}$;
\item[\rm{(3)}] on the singular fibre $\cS_0$, $\Omega|_{X^+}=\omega^+$,\  \  $\Omega|_{X^-}=\omega^-$;
\item[\rm{(4)}] $J_\cS$ preserves $T\cS_\lambda$ for every $\lambda\in D^*$;
\end{itemize}
\et
Then every $(\cS_\lambda,\omega_\lambda)$, $\lambda\in D^*$, is a smooth compact symplectic manifold, 
called symplectic sum of $X^+$ and $X^-$ along $Z$ with gluing parameter $\lambda$, and simply denoted as $\cS_\lambda=X^+\#_{Z,\lambda} X^-$. 
They are symplectically isotopic to one another and can be regarded  as deformations, in the symplectic category, of the singular fiber $\cS_0$.

\subsection{Normal form for  a codimensional two symplectic submanfiold}

Let $Z$ be  a smooth  manifold  and $\pi_N:(N,\i_N)\to Z$  be a complex line 
bundle (a rank-2 real vector bundle with a complex structure $\i_N$ on each fiber).
Let $(g_N,\nabla^N)$ be a  Hermitian structure on $(N,\i_N)$, which consists of  a Hermitian metric and a  connection on $N$
such that for $\forall v,w\in N_z$, $z\in Z$  and any  sections $\xi,\eta\in\Gamma(N)$, the following hold
\bq\label{Her1}
g_N(\i_N v,w)=\i g_N( v,w)=-g_N(v,\i_Nw),
\eq
\bq\label{Her2}
\nabla^N (\i_N\xi) =\i_N\nabla^N \xi,
\eq
\bq\label{Her3}
d\big(g_N(\xi,\eta)\big)=g_N(\nabla^N\xi,\eta)+g_N(\xi,\nabla^N\eta).
\eq
Denote by 
$$\rho_N: N\to\R,$$
\bq\label{5def-rho}
\rho_N ( v ) = g_N ( v , v ) = | v |^2
\eq 
the square of the norm function. Let 
$$\pi_{SN} : SN\to Z$$ be the circle bundle of $N$.
The connection $\nabla^N$ induces a splitting of  the following short  exact sequence
\bq\label{5N-exact}
0 \longrightarrow  T^{vrt}(SN)\cong \ker(d\pi_{SN})\longrightarrow T(SN)\xlongrightarrow{d\pi_{SN}}  \pi_{SN}^*  TZ\longrightarrow 0
\eq
of vector bundles over $SN$. Denote by $\alpha_{SN}$ the  1-form on $SN$ vanishing on
the image of $\pi_{SN}^*TZ$ in $T(SN)$ corresponding to this splitting such that
\bq\label{4alp-SN}
\alpha_{SN} \left(\frac{d}{d\theta} e^{\i\theta}v \Big|_{\theta=0}\right)   =1,\ \ \  \forall\ v\in SN.
\eq
It can be extended  to a  1-form $\alpha_N$ on $N - Z$ via the radial retraction
$$N- Z\to SN, \ \ \ \   v\to \frac{v}{|v|}.$$
Then the 1-form $\rho_N \alpha_N$ is  well-defined and smooth on the total space of the line bundle $N$.

Assume now that $(Z,\omega_Z)$ is a symplectic manifold.
For $\epsilon>0$,  define  a  2-form on the total space of a complex line bundle $N$ over $Z$ 
\bq\label{5form-Z}
\omega^\epsilon_{N,Z}:=\pi_N^*\omega_Z+\frac{\epsilon^2}{2}d(\rho_N \alpha_N).
\eq
We see that $\omega^\epsilon_{N,Z}$ is closed and the restriction to the zero section of $N$:  $\omega^\epsilon_{N,Z}|_{Z}=\omega_Z$.
If $Z$ is compact, since $d\rho_N\wedge\alpha_N$ is nondegenerate on each fiber, there exists $\epsilon_Z\in\R^+$ such 
that the restriction of $\omega^\epsilon_{N,Z}$ to  
$N(\delta)=\{ v\in N:\ |v|< \delta \}$ is  non-degenerate whenever $\delta,\epsilon\in\R^+$ and $\delta\epsilon<\epsilon_Z$.
So $(N(\delta) , \omega^\epsilon_{N,Z})$ is a symplectic manifold  if $\delta\epsilon <\epsilon_Z$.  
In particular, if $\delta<\epsilon_Z$, then 
$(N(\delta) , \omega^1_{N,Z})$ is a symplectic manifold.

Now consider a symplectic manifold $(X,\omega_X)$. Let $Z$ be a codimension 2 symplectic submanifold of $X$.
Denote the symplectic normal bundle of $Z$ in $X$ by
\bq\label{5iso-NZ}
N_XZ:=\frac{TX|_Z}{TZ}\cong (TZ)^{\omega_X},
\eq where 
$$ TZ^{\omega_X}=\{v\in T_xX\ | \  x\in Z, \omega_X(v,w)=0,\ \forall\ w\in T_xZ\}.$$
So under the identification (\ref{5iso-NZ}) we have a symplectic orthogonal decomposition 
$$
TX|_Z=TZ\oplus N_XZ,\ \ \ \  \omega_X=\omega_Z\oplus\omega_N,
$$
where $\omega_N$ is a fiberwise symplectic structure on $N_XZ$ such that on each fiber $N_z$, $\omega_N=\omega_X|_{N_z}$.
Let $\omega_X^N$ be the induced symplectic form on the  normal bundle $N_XZ$.  
We can  especially consider (fiberwise) {\bf $\bm\omega_X^N$-compatible} complex structures $\i_X$ on $N_XZ$  
and  {\bf$\bm\omega_X^N$-compatible} Hermitian structures $(g_X,\nabla^X)$ on  $(N_XZ , \i_X )$, $i.e.$
\bq\label{Her-comp}
\omega_X^N(\i_Xv,\i_Xw)=\omega_X^N(v,w),\ \ \ \   g_X(v,w)=\omega_X^N(v,\i_Xw)
\eq 
for $\forall\ v,w\in N_XZ|_z, z\in Z$.
The spaces of (fiberwise) $\omega_X^N$-compatible complex structures on $N_{X}Z$  and of $\omega_X^N$-compatible Hermitian structures on 
$(N_{X}Z,\i_X)$ are non-empty and contractible.

For example, when $\delta\epsilon<\epsilon_N$ is sufficiently small and  as a zero section, $Z$ is a compact codimension 2 symplectic submanifold of  
$(N(\delta) , \omega^\epsilon_{N,Z})$,   the symplectic form  $\omega^\epsilon_{N,Z}$ in (\ref{5form-Z}) 
is just the form $\omega_X^N|_{N(\delta)}$ up to a symplectomorphism.
Then we have
$$TZ^{\omega^{\epsilon}_{N,Z}}\cong N_{N(\delta)}Z,$$
$$(\omega^{\epsilon}_{N,Z})_N=\frac{\epsilon^2}{2}d(\rho_N \alpha_N)|_{N(\delta)}.$$
 
Note that when $\delta<\epsilon_Z$ is sufficiently small, the 2-form $\omega^1_{N,Z}$ restricts to $\omega^N_X$  on $T(N_{Z,X})|_Z$ under the isomorphism as in (\ref{5iso-NZ}).
Then by  the Symplectic Neighborhood Theorem, 
a neighborhood of $Z$ is completely determined by the restriction of $\omega$ to $Z$ together with 
the isomorphism class of the symplectic normal bundle $N_{X}Z$.
That is, there exist $\delta_Z>0$, an smooth injective open map,
\bq\label{local:diff}
\Psi_X: (N_{X}Z(\delta_Z),Z)\longrightarrow (X, Z)
\eq
 such that
$$d_x\Psi_X=id, \ \forall\ x\in Z,\ \ \ \   \Psi_X^*\omega_X=\omega^1_{N,Z}|_{N_{X}Z(\delta_Z)}.$$  Note that here  we can choose 
 $\delta_Z < \epsilon_Z$ so that $\omega^1_{N,Z}|_{N_{X}Z(\delta_Z)}$ is a symplectic form on $N_{X}Z(\delta_Z)$.

For any $\epsilon>0$, define 
\bq\label{5PXE}
\Psi_{X;\epsilon}: (N_{X}Z(\epsilon^{-1}\delta_Z),Z)\longrightarrow (X,Z)
\eq 
$$\Psi_{X;\epsilon}(p=(z,v))=\Psi_X(z,\epsilon v),\ \ \  z\in Z, \ v\in N_{X}Z(\epsilon^{-1}\delta_Z)|_z,$$ then  $\Psi_{X;\epsilon}$  is a   smooth injective open map,  satisfies 
$$\Psi_{X;\epsilon}^*\omega_X=\omega^\epsilon_{N,Z}|_{N_{X}Z(\epsilon^{-1}\delta_Z)}$$
and restricts to the identity on  $Z$.

\subsection{Symplectic  sum }\label{G3p}

Let $(X^+, \omega^+)$ and $(X^-, \omega^-)$ be two $2n$-dimensional compact symplectic manifolds and $Z\subset X^+, X^-$ 
be a common symplectic codimension 2 submanifold satisfying  $\omega_Z=\omega^+|_Z=\omega^-|_Z$ so that (\ref{5EL}) holds.
Fix (fiberwise) complex structures $\i_{+}$ and $\i_{-}$ on the normal bundles
$$\pi_{+,Z} :N_{+}Z\to Z\ \ \  \mathrm{and}\ \ \ \pi_{-,Z} :N_{-}Z\to Z$$
that are compatible with $\omega_+^N$ and $\omega_-^N$, respectively. 
Here $\omega_\pm^N=\omega_{X^\pm}^{N_\pm Z}$ as in the last subsection.
Fix  an isomorphism $\Phi$ of trivial complex line bundles (\ref{5iso-LB}),  let  
$$\Phi_2:N_{+} Z \otimes_\C  N_{-} Z\overset{\Phi}{\cong}Z\times\C \to \C$$ 
be the composition of $\Phi$ with the projection $Z\times\C\to\C$.
\bdf\label{phi-compat}
We say that an isomorphism $\Phi$ is $(\omega_+^N, \omega_-^N)$-{\bf compatible} if
\bq\label{phi2}
 |\Phi_2\big( (z,v)\otimes (z,w) \big)|^2 = |\Phi_2( v\otimes  w)|^2 =\omega_+^N(v,\i_+v)\cdot \omega_-^N(w,\i_-w),
\eq
for $\forall\ v\in N_{+}Z|_z$,  $w\in N_{-}Z|_z$,  $z\in Z$.
\edf
We choose an isomorphism $\Phi$ which is  $(\omega_+^N, \omega_-^N)$-compatible.  
In fact, (\ref{phi2})  can be achieved by scaling any given isomorphism $\Phi$ in  (\ref{5iso-LB}) and this does not
change the homotopy class of $\Phi$.

Choose Hermitian structures $(g_+ , \nabla_+ )$ on $(N_{+}Z, \i_+ )$ and $(g_- , \nabla_- )$ on $(N_{-}Z, \i_- )$ 
that are compatible with $\omega_+^N$ and $\omega_-^N$ , in the sense of (\ref{Her-comp}), and compatible with $\Phi_2$, in the following sense 
\bq\label{H-compat}
 |\Phi_2(z, v\otimes  w)|^2 = \rho_{+}(v)\cdot \rho_{-}(w),
\eq
\bq\label{H-compat-d} 
d\big(\Phi_2(\xi\otimes_\C \eta)\big) = \Phi_2\big((\nabla_+\xi) \otimes_\C \eta\big) + \Phi_2\big(\xi\otimes_\C (\nabla_-\eta)\big) 
\eq
for $\forall\ v\in N_{+}Z|_z$,  $w\in N_{-}Z|_z$,  $z\in Z$ and $\xi\in\Gamma(Z;N_{+}Z)$, $\eta\in\Gamma(Z;N_{-}Z)$.
    Here    $\rho_\pm$  is the fiberwise normal square functions on $N_\pm Z$, and the equation (\ref{H-compat-d} ) is an identity of
    differential 1-forms on $Z$. 
 
 Denote by $\alpha_+$ and $\alpha_-$ the connection 1-forms on $N_{+} Z-Z$ and $N_{-} Z-Z$ corresponding to
$(g_+ , \nabla_+ )$  and $(g_- , \nabla_- )$, respectively. For a sufficiently small  $\epsilon>0$, as  in (\ref{5form-Z}), define
\bq\label{form-Z-pm}
\omega^\epsilon_{+,Z}:=\pi_{+,Z}^*\omega_Z+\frac{\epsilon^2}{2}d(\rho_+ \alpha_+),\  \ \ \   \ 
\omega^\epsilon_{-,Z}:=\pi_{-,Z}^*\omega_Z+\frac{\epsilon^2}{2}d(\rho_- \alpha_-).
\eq
 As in (\ref{local:diff}),  there exist    $\delta_Z>0$  and the smooth injective  open maps 
\[
\Psi_\pm: (N_{\pm}Z(\delta_Z), Z)\longrightarrow (X^\pm, Z)
\]
such that 
$$d_x\Psi_\pm=id, \ \forall\ x\in Z,\ \ \ \   \Psi_\pm^*\omega_{X^\pm}=\omega^1_{\pm,Z}|_{N_{\pm}Z(\delta_Z)}.$$ 
 Following (\ref{5PXE}),  there exist  smooth injective open maps
 \bq\label{PXE-pm}
\Psi_{+;\epsilon}: (N_{+}Z(\epsilon^{-1}\delta_Z),Z)\longrightarrow (X^+,Z), \ \ \ \  \Psi_{-;\epsilon}: (N_{-}Z(\epsilon^{-1}\delta_Z),Z)\longrightarrow (X^-, Z).  
\eq 
  satisfying 
\bq\label{Pxi-symp}
\Psi_{+;\epsilon}^*\omega_+=\omega^\epsilon_{+,Z}|_{N_{+}Z(\epsilon^{-1}\delta_Z)},\ \ \ \   
\Psi_{-;\epsilon}^*\omega_-=\omega^\epsilon_{-,Z}|_{N_{-}Z(\epsilon^{-1}\delta_Z)}
\eq
and restricting  to the identity on  $Z$.     Note that   $\epsilon$ and $\delta$ are chosen, $\epsilon^{-1}\delta_Z >2$, so that 
$N_{\pm}Z(2)$  is contained in   the domain
of $\Psi_{\pm;\epsilon}$.

Consider  the projections
$$\Pi_Z:N_{+}Z\oplus N_{-}Z \longrightarrow Z,$$
$$\Pi_+:N_{+}Z\oplus N_{-}Z \longrightarrow N_{+}Z,$$
$$\Pi_-:N_{+}Z\oplus N_{-}Z \longrightarrow N_{-}Z, $$
and the natural product map
\ba\label{prod-map}
\cP:N_{+}Z\oplus N_{-}Z &\to& \C,\nonumber\\
(z, v, w) &\mapsto& \Phi_2(v\otimes w).
\na
Then the  following identity  holds\footnote{This identity is just the formula (3.8) of \cite{FZ}.}
\bq\label{P-omega}
\cP^* \omega_\C=\frac{1}{2}d(\rho_+\rho_-\cP^*d\theta)=\frac{1}{2}d\big(\rho_+\rho_-(\Pi_+^*\alpha_+ + \Pi_-^*\alpha_-)\big),
\eq
where $\omega_\C=rdr\wedge d\theta=\frac{1}{2}d(r^2d\theta)$ is the standard symplectic form on $\C$, 
$\rho_+$ and $\rho_-$ also denote the  extensions on $N_{X^+}Z\oplus N_{X^-}Z$ defined by $\rho_+(v,w)=|v|^2$,  $\rho_-(v,w)=|w|^2$.

The following  pieces are the basic  building blocks in the symplectic sum construction
\ba\label{E-pm}
\cE_+ &:=& \big(X^+ - \Psi_{+;\epsilon}(\overline{N_{+}Z(1)}) \big)\times \D_{\delta}, \\
\cE_- &:=& \big(X^- - \Psi_{-;\epsilon}(\overline{N_{-}Z(1)}) \big)\times \D_{\delta} ,  \\
 \cS_Z&:=& \{(z,v,w)\in N_{+}Z\oplus N_{-}Z\ \big|\ |v|,|w|< 2, \epsilon|\cP(v,w)|<\delta\},\\
 \cS_{Z,+}&:=&\{(z,v,w)\in \cS_Z: |v|>1\},\ \ \   \ \ \     \cS_{Z,-}:=\{(z,v,w)\in \cS_Z: |w|>1\}.
\na
where $ \D_{\delta}=\{\zeta\in\C:\ |\zeta|<\delta\}$, and $\overline{N_{\pm}Z(1)}$ denote the closed unit disc bundles of $N_{\pm}Z$.  
  We first choose  a sufficiently small $\epsilon>0$, then assume that 
\bq\label{E-D}
2\epsilon<\delta_Z,\ \ \ \    2\delta<\epsilon.
\eq
 Define the gluing maps to be the open maps
 \bq\label{gl+}
gl_+: \cS_{Z,+}\longrightarrow \cE_+, \ \ \ \   (z,v,w)\mapsto\big(  \Psi_{+;\epsilon}(z,v), \epsilon\cP(v,w) \big), 
\eq
\bq\label{gl-}
gl_-:  \cS_{Z,-}\longrightarrow  \cE_-, \ \ \ \   (z,v,w)\mapsto\big(  \Psi_{-;\epsilon}(z,w), \epsilon\cP(v,w) \big).
\eq
From the assumptions in (\ref{E-D}), $gl_+$ and $gl_-$ are well-defined diffeomorphisms between open subsets of their domains and targets. 
Let $\cS$ be the resulting smooth manifold from gluing $\cE^+$, $\cE^-$ and $\cS_Z$ by the  maps $gl_\pm$
\bq\label{S}
\cS:=\cE_{-}\ \bigcup_{gl_-}\ \cS_Z\ \bigcup_{gl_+}\ \cE_{+}. 
\eq 
The maps
\[
\pi_{\pm,\C}: \cE_\pm  \longrightarrow \D_\delta,\ \ \ : (x^\pm,\lambda)  \mapsto \lambda,\]
\[
\pi_{Z,\C}: \cS_Z \longrightarrow \D_{\delta}, \ \ \  (z,v,w)\mapsto \epsilon\cP(v,w)
\]
are intertwined by $gl_+$ and $gl_-$, so they induce a smooth map
\bq\label{Piep}
\pi_\epsilon:\cS\longrightarrow \D_{\delta}.
\eq
By the second assumption in (\ref{E-D}), every fiber $\cS_\lambda=\pi_\epsilon^{-1}(\lambda)$ of $\pi_\epsilon$  
  is compact.

 We now follow \cite{FZ} to  construct a symplectic form $\Omega_\cS^\epsilon$ on $\cS$.  
First,   note that  $\cE_\pm  $ are symplectic with the symplectic form given by  $ p_\pm^*\omega_+  +  \pi_{+,\C}^*\omega_\C$ 
where $p_\pm$ is the natural projection  
\bq  \label{pi-pm-X}
p_\pm: \cE_\pm\longrightarrow  X^\pm - \Psi_{\pm;\epsilon}(\overline{N_{X^\pm}Z(1)}),\ \ \ \  :  (x^\pm,\lambda) \mapsto x^\pm. 
\eq
Take a cut-off function $\beta:\R\to[0,1]$, which is a smooth function such that
\bq\label{beta}
\beta(t)=\begin{cases} 0, & \mathrm{if} \  t\le\frac12 ,\\    1, & \mathrm{if} \  t\ge 1.\end{cases}
\eq
 Define a closed 2-form on $\cS_Z$
\ba\label{omgZ}
\omega^\epsilon_{\cS_Z} :=  \Pi_{Z}^*\omega_Z+\frac{\epsilon^2}{2} d (\alpha_\#),
\na
where the 1-form
\ba\label{alp-SZ}
\alpha_\# &:= &  (1-\beta\circ\rho_-)\Pi_+^*(\rho_+ \alpha_+) + (1-\beta\circ\rho_+)\Pi_-^*(\rho_- \alpha_-)\nonumber\\
 &  & \ \ \    + (\beta\circ\rho_+ + \beta\circ\rho_-)\rho_+\rho_-( \Pi_+^*\alpha_+ + \Pi_-^*\alpha_- ).
\na

Since  on $\cS_Z$, $|\cP(v,w)|<\epsilon^{-1}\delta<\frac12$, we have  $\rho_+^{\frac12}=|v|>1$, $\rho_-^{\frac12}=|w|<\frac12$ on $\cS_{Z,+}$.  In particular, $\beta\circ \rho_+ =1$ and $\beta\circ \rho_- =0$ on $\cS_{Z,+}$.
Using (\ref{Pxi-symp}), (\ref{P-omega}) and (\ref{gl+}),  the restriction of this closed 2-form $\omega^\epsilon_{\cS_Z}$ to $\cS_{Z,+}$ is
\ba\label{res=gl+}
\Pi_{Z}^*\omega_Z&+&\frac{\epsilon^2}{2}d \Big( \Pi_+^*(\rho_+ \alpha_+)+ \rho_+\rho_-( \Pi_+^*\alpha_+ + \Pi_-^*\alpha_- )\Big)\nonumber\\
&=&\big( \Pi_{Z}^*\omega_Z+\frac{\epsilon^2}{2}d  \Pi_+^*(\rho_+ \alpha_+)\big) +\epsilon^2\frac{1}{2} d \big( \rho_+\rho_-( \Pi_+^*\alpha_+ + \Pi_-^*\alpha_- )\big)\nonumber\\
&=& gl_+^* (p^*_+\omega_+)  +  \epsilon^2\cP^*\omega_\C= gl_+^* (p^*_+ \omega_+  +  \pi_{+,\C}^*\ \omega_\C). 
\na
Similarly, the restriction of $\omega^\epsilon_{\cS_Z}$ to $\cS_{Z,-}$ is  $ gl_-^* (p^*_- \omega_-  +  \pi_{-,\C}^*\ \omega_\C)$. 
So along with the 2-forms 
$$p^*_+\omega_+  +  \pi_{+,\C}^*\ \omega_\C\ \mathrm{on}\ \  \cE_+\ \ \ \mathrm{and}\ \ \    p^*_-\omega_-  + \pi_{-,\C}^*\ \omega_\C\  \mathrm{on}\ \  \cE_-,$$ 
we get a closed 2-form 
$\Omega_\cS^\epsilon$ on $\cS$ via the  gluing construction (\ref{S}).

From the same calculation as the one of \cite[Page 23]{FZ}, one can verify that the closed 2-form (\ref{omgZ}) 
is also nondegenerate if  $\epsilon$  and  $\delta$  is chosen to be small enough as in (\ref{E-D}).
Thus, we obtained a symplectic manifold $(\cS,\Omega=\Omega^\epsilon_{\cS})$ of dimension $2n+2$ such that
\bq\label{OMG}
\Omega^\epsilon_{\cS}=
\begin{cases} p^*_-\omega_-  +  \pi_{-,\C}^*\ \omega_\C, & \mathrm{on} \  \cE_-, \\  \ \ \ \ \ \   \omega_{\cS_Z}^\epsilon, & \mathrm{on} \  \cS_Z, 
\\  p^*_+\omega_+  +  \pi_{+,\C}^*\ \omega_\C, & \mathrm{on} \   \cE_+,  
\end{cases}
\eq
and a fibration (see (\ref{Piep}))
\bq\label{Pi-sum}
\pi := \pi_\epsilon:   \cS\longrightarrow \D_{\delta}.
\eq

\bigskip

Note that 
$$gl_+^{-1}\Big[ \Big(\Psi_{+;\epsilon}\big(N_{+}Z(2) - \overline{N_{+}Z(1)} \big) \Big)\times\{0\}\Big]=\{(v,0)\in N_{+}Z\oplus \{0\}\ \big|\ 1<|v|< 2\},$$
$$gl_-^{-1}\Big[ \Big(\Psi_{-;\epsilon}\big(N_{-}Z(2)- \overline{N_{-}Z(1)} \big) \Big) \times\{0\}\Big]=\{(0,w)\in \{0\}\oplus N_{-}Z\ \big|\ 1<|w|< 2\},$$
and
$$\pi_{Z,\C}^{-1}(0)=\{(v,0)\in N_{+}Z\oplus \{0\}\ \big|\ |v|< 2 \}\ \bigcup\ \{(0,w)\in \{0\} \oplus N_{-}Z\ \big|\ |w|< 2 \}.$$
So 
\ba
\cS_0=\pi^{-1}(0) &:= & \big(X^+ - \Psi_{+;\epsilon}(\overline{N_{+}Z(1)}) \big) \times \{ 0\}\ \bigcup_{gl_+}\  \pi_{Z,\C}^{-1}(0) \nonumber\\ 
                            &     & \bigcup_{gl_-} \big(X^- - \Psi_{-;\epsilon}(\overline{N_{-}Z(1)}) \big)\times \{0\} \nonumber\\
                      & \cong & X^+\  \bigcup_Z\  X^-
\na 
is a singular fiber of (\ref{Pi-sum}).
For $\lambda\in \D_{\delta}^*=\D_{\delta}\setminus\{0\}$,  denote the fibre by $\cS_\lambda=\pi^{-1}(\lambda)$. To study the  nondegeneracy of the restriction of $\Omega^\epsilon_{\cS}$ to the fibre $\cS_\lambda$ 
(the restriction of $\omega^\epsilon_{\cS_Z}$ to the fibre $\cS_\lambda\cap\cS_Z$), 
one needs to construct $\Omega^\epsilon_{\cS}$-tame and compatible almost complex structure $J_{\cS}$ on $\cS$ which 
preserves the tangent spaces to the fibers of the fibration (\ref{Pi-sum}). 
We refer the reader to \cite{FZ} for the detailed construction and the proof  the following proposition. 
\begin{proposition}\label{Prop-sub-J}
For every $ \lambda\in \D^*_{\delta}$ with sufficiently small $\delta$,  the fibre
$\cS_\lambda $ is a compact symplectic submanifold of $(\cS,\Omega^{\epsilon}_{\cS})$ with restricted symplectic form 
$\omega_\lambda:=\Omega^{\epsilon}_{\cS}|_{\cS_\lambda}$, 
called the symplectic sum of $X^+$ and $X^-$ along $Z$, and symplectically isotopic to one another. 
Moreover,  there exists a compatible almost complex structure  $J_{\cS}$ on $\cS$, such that  $\cS_\lambda$ is $J_\cS$-invariant. 
\end{proposition}

\section{Sum of Hamiltonian manifolds}\label{SEC-Sum-Ham}

 In this section,  we establish the main  result of the paper on the operation of sum of two Hamiltonian manifolds along a common
codimension 2 Hamiltonian submanifold.  
Fix a compact and connected Lie group $G$.  Consider a Hamiltonian G-manifold $(X,\omega_X,\mu_X)$.
Let $Z$ be a codimension 2 compact  Hamiltonian submanifold of $X$.   Under the identification (\ref{5iso-NZ}) we have a  $G$-equivariant symplectic orthogonal decomposition 
$$
TX|_Z=TZ\oplus N_XZ,\ \ \ \  \omega_X=\omega_Z\oplus\omega_N, 
$$
where $N_XZ = (TZ)^{\omega_X}$ is the $G$-equivariant  symplectic normal bundle of $Z$ in $X$.  

The induced symplectic form   $\omega_X^N$  on the total space of the normal bundle $\pi_N: N_XZ \to Z$  is $G$-invariant.    The spaces of (fiberwise) $\omega_X^N$-compatible complex structures on $N_{X}Z$  and of $\omega_X^N$-compatible Hermitian  metrics  on 
$(N_{X}Z,\i_X)$ are non-empty and contractible. By averaging over $G$, we can choose a $G$-invariant 
$\omega_X^N$-compatible complex structure   $\i_X$ on $N_{X}Z$,  a $G$-invariant  $\omega_X^N$-compatible Hermitian  metric  $g_N$   and $G$-invariant Hermitian connection $\nabla^N$  on   $(N_{X}Z,\i_X)$.  
In particular,  the square of the norm function  
\[
\rho_N: N_X Z \to\R, \qquad \rho_N (x, v) = g_N(v, v)  \ \ \forall\ v \in N_XZ|_x, x\in Z
\]
 is $G$-invariant.  The  $G$-invariant connection  $\nabla^N$  on $N_{X}Z$  induces a $G$-invariant  connection 1-form on the unit circle bundle  $S N_X$ of $ N_X Z$,  and can be extended to a  $G$-invariant 1-form $\alpha$ on $N_X Z  -Z$ via the radial retraction 
\[
N_X Z  -Z \longrightarrow S N_X, \qquad v\mapsto \dfrac{v}{|v|}. 
\]
Therefore  $\rho_N \alpha_N$ is  a $G$-invariant 1-form on the total space of $N_X Z$. 
If  $\epsilon\delta$ is small enough, then  the restriction of 
\[
\omega^\epsilon_{N,Z} :=\pi_N^*\omega_Z+\frac{\epsilon^2}{2}d(\rho_N \alpha_N).
\]
to  $N_{X}Z(\delta )  =\{ v\in N:\ |v|< \delta\}$ is  a $G$-invariant symplectic form.
From the discussion in  subsection \ref{OA-mmp}, we know that the $G$-action on $N_{X}Z(\delta )$ is Hamiltonian with the moment map 
$\mu_{N, Z}^\epsilon:  N_{X}Z(\delta ) \to \g^*$ 
given by 
\[
\mu_{N, Z}^\epsilon = \mu_Z\circ \pi_N + \frac{\epsilon^2}{2} \mu_{\rho_N \alpha_N},
\]
where $\mu_{\rho_N \alpha_N}$ is the $\rho_N \alpha_N$-moment map   defined in subsection \ref{OA-mmp}. 

Since $Z$ is $G$-invariant in both $N_{X}Z(\delta )$ and $X$, and $\omega^\epsilon_{N,Z}|_Z=\omega_X|_Z$, 
by an equivariant version of the Symplectic Neighborhood Theorem(a direct result from the Darboux-Weinstein Theorem in \cite{GS}),  
there exist $\delta_Z>0$, a tubular neighborhood $\cO_X(Z)$ of $Z$ in $X$ and a $G$-equivariant  diffeomorphism
$$\Psi_X: (N_{X}Z(\delta_Z), Z)\longrightarrow (\cO_X(Z), Z)$$ such that
$$   \Psi_X^*\omega_X=\omega^1_{N,Z}|_{N_{X}Z(\delta_Z)},\ \ \   \Psi_X|_Z = Id_Z.$$  Moreover, we have
\[
\mu_X \circ  \Psi_X =   \mu_Z\circ \pi_N + \frac{1}{2} \mu_{\rho_N \alpha_N}. 
\]
For any $\epsilon>0$,  the  smooth  $G$-equivariant injective open map \bq\label{PXE}
\Psi_{X;\epsilon}: (N_{X}Z(\epsilon^{-1}\delta_Z),Z)\longrightarrow (X,Z)
\eq 
$$\Psi_{X;\epsilon}(z,v)=\Psi_X(z,\epsilon v),\ \ \  z\in Z, \ v\in N_{X}Z(\epsilon^{-1}\delta_Z)|_z,$$
satisfies 
$$\Psi_{X;\epsilon}^*\omega_X=\omega^\epsilon_{N,Z}|_{N_{X}Z(\epsilon^{-1}\delta_Z)}$$
and restricts to the identity on  $Z$.  The corresponding moment maps 
satisfy
\bq\label{mu-dZ}
\mu_X \circ  \Psi_{X;\epsilon}=\mu_Z\circ\pi_{N}+ \frac{\epsilon^2}{2} \mu_{\rho_N\alpha_N}. 
\eq

 Let $(X^+,\omega^+,\mu^+)$ and $(X^-,\omega^-,\mu^-)$ be two Hamiltonian $G$-manifolds of the  dimension $2n$, 
$Z=Z^+ \cong Z^-$ be their common codimension 2 Hamiltonian $G$-submanifold satisfying 
\[
\omega_Z = \omega^+|_Z = \omega^-|_Z, \qquad \mu_Z=\mu^+|_Z=\mu^-|_Z.
\]
 Take $G$-invariant  almost complex structures  $J^\pm$ on $X^\pm$ compatible with $\omega^\pm$ such that
 $J^+|_{TZ} = J^-|_{TZ}$.  The symplectic normal bundles 
 \[
 \pi_{\pm , Z}:  N_\pm Z = (TZ)^{\omega^\pm}  \longrightarrow Z
 \]
are  $G$-equivariant  and  can be endowed with $G$-invariant fiberwise complex structures $\i_\pm$, compatible with the induced 
$G$-invariant symplectic form $\omega_\pm^N = \omega^\pm|_{ N_\pm Z}$ as above.  
In particular,   these determine $G$-invariant Hermitian metrics $g_\pm$ on $(N_\pm Z, \i_\pm)$ in the sense of (\ref{Her-comp}).   
We assume that $N_+ Z$ and $N_- Z$ have opposite $G$-equivariant Euler classes:
\bq\label{EL}
e_G(N_{+}Z)+e_G(N_{-}Z)=0.
\eq
By the classification theorem of equivariant complex line bundles (Theorem C.47 \cite{GGK}), there is a 
$G$-equivariant  isomorphism  of complex line bundles
\bq\label{iso-LB}
\Phi:N_{+} Z \otimes_\C N_{-} Z \longrightarrow Z\times\C,
\eq
where $Z\times \C$ is a trivial $G$-equivariant line bundle in the sense that $G$-action on $\C$ is trivial. 
  
 Fix a $G$-equivariant isomorphism $\Phi$ 
which is compatible with $(\omega_+^N, \omega_-^N)$ (see (\ref{phi2})).  
Choose   $G$-invariant  $\omega_\pm^N$-compatible connections $\nabla_\pm$ on $N_\pm Z$ in the sense of (\ref{H-compat}) and (\ref{H-compat-d}).  Then we have $G$-invariant norm square functions
\[
\rho_\pm:  N_\pm Z \longrightarrow \R, \qquad  \rho_\pm (v) = g_\pm (v, v). 
\]
Choose $G$-invariant  connections $\nabla_\pm$ on $N_\pm Z$ which in turn  define $G$-invariant connection 1-forms $\alpha_\pm$ on the unit circle $SN_\pm Z$ of $N_\pm Z$.  Note that $\alpha_\pm$ can  be extended to   $G$-invariant 1-forms
on $N_\pm Z-Z$ via the radial retraction. Therefore we have 
 $G$-invariant 1-forms  $\rho_\pm \alpha_\pm$ on the total spaces of the symplectic normal bundles  $N_\pm Z \to Z$.   
 
 When $\epsilon \delta$ is small enough,
\[
\big(N_\pm Z (\delta), \omega_{\pm, Z}^\epsilon = \pi_{\pm , Z}^* \omega_Z + \frac{\epsilon^2}{2} d (\rho_\pm\alpha_\pm) \big) 
\]
are  Hamiltonian $G$-manifolds with the moment map $\mu_{\pm, Z}^\epsilon:  N_\pm Z (\delta) \to \g^*$
given by
\[
\mu_{\pm, Z}^\epsilon = \mu_Z\circ \pi_{\pm, Z} + \frac{\epsilon^2}{2} \mu_{\rho_\pm  \alpha_\pm},
\]
where $\mu_{\rho_\pm  \alpha_\pm}$  are the $\rho_\pm \alpha_\pm$-moment maps   defined in subsection \ref{OA-mmp}. 

 There exist $\delta_Z>0$,   tubular neighborhoods  $\cO_{X^\pm} (Z)$ of $Z$ in $X^\pm$ and $G$-equivariant  diffeomorphisms 
$$\Psi_\pm : (N_{\pm}Z(\delta_Z), Z)\longrightarrow (\cO_{X^\pm} (Z), Z)$$ such that
$$   \Psi_\pm ^*\omega^\pm =\omega^1_{\pm ,Z}|_{N_{\pm }Z(\delta_Z)},\ \ \   \Psi_\pm |_Z = Id_Z.$$  Moreover, we have
\[
\mu_X \circ  \Psi_X =   \mu_Z\circ \pi_N + \frac{1}{2} \mu_{\rho_N \alpha_N}. 
\]
  and   for any $\epsilon>0$ the  smooth  $G$-equivariant injective open maps of Hamiltonian $G$-manifolds
 \[
 \Psi_{\pm;\epsilon}:  N_{\pm}Z(\epsilon^{-1}\delta_Z)  \longrightarrow   X^\pm, 
\]
$$\Psi_{\pm;\epsilon}(z,v)=\Psi_\pm (z,\epsilon v),\ \ \  z\in Z, \ v\in N_{\pm}Z(\epsilon^{-1}\delta_Z)|_z,$$
satisfying 
$$\Psi_{\pm ;\epsilon}^*\omega_\pm=\omega^\epsilon_{\pm, Z}|_{N_{\pm }Z(\epsilon^{-1}\delta_Z)}$$
and restricts to the identity on  $Z$.  The corresponding moment maps 
satisfy
\bq
\mu_X \circ  \Psi_{X;\epsilon}=\mu_Z\circ\pi_{N}+ \frac{\epsilon^2}{2} \mu_{\rho_N\alpha_N}. 
\eq

Recall the symplectic sum construction (\ref{S})
\[
\cS=\cE_{-}\ \bigcup_{gl_-}\ \cS_Z\ \bigcup_{gl_+}\ \cE_{+}. 
\]
Here  $\cE_\pm   =  \big(X^\pm  - \Psi_{\pm ;\epsilon}(\overline{N_{\pm }Z(1)}) \big)\times \D_{\delta}$ with symplectic form $p^*_\pm \omega^\pm   +  \pi_{\pm ,\C}^*\ \omega_\C$.   Note that $G$-action on $\D_\delta$ is trivial, 
  the  symplectic form $p^*_\pm \omega^\pm   +  \pi_{\pm ,\C}^*\ \omega_\C$ is $G$-invariant. Therefore,    $( \cE_\pm,   p^*_\pm \omega^\pm   +  \pi_{\pm ,\C}^*\ \omega_\C)$ are Hamiltonian $G$-manifolds
with   the moment maps $\mu_{\cE_\pm}:  \cE_\pm  \to \g^*$  given by 
\[
\mu_{\cE_\pm} =  \mu^{\pm} \circ p_\pm.  
\]
The  symplectic form  $\omega^\epsilon_{\cS_Z} =    \Pi_{Z}^*\omega_Z+\frac{\epsilon^2}{2} d (\alpha_\#)$  on the middle part $ \cS_Z =    \{(z,v,w)\in N_{+}Z\oplus N_{-}Z\ \big|\ |v|,|w|< 2, \epsilon|\cP(v,w)|<\delta\}$  
 is $G$-inavriant as  $\alpha_\#$,   defined in  (\ref{alp-SZ}).   This implies that 
 \[ ( \cS_Z, \omega^\epsilon_{\cS_Z})
 \]
 is a  Hamiltonian $G$-manifold with 
    the associated moment map   $\mu_{\cS_Z } :  \cS_Z \to \g^*$  given by 
\[
\mu_{\cS_Z } = \mu_Z \circ   \Pi_{Z} + \frac{\epsilon^2}{2} \mu_{\alpha_\#}.
\]
Under the gluing maps  $gl_\pm:  \cS_{Z, \pm} \to \cE_\pm$  defined in  (\ref{gl+}) and (\ref{gl-}), we have 
\[
gl_\pm^* (p^*_\pm  \omega^\pm   +  \pi_{\pm,\C}^*\ \omega_\C) = \omega^\epsilon_{\cS_Z} |_{ \cS_{Z, \pm}}, 
\]
which implies 
\[
\mu_{\cS_Z }|_{ \cS_{Z, \pm} }= \mu_{\cE_\pm} \circ gl_\pm =  \mu^{\pm} \circ p_\pm \circ gl_\pm.
\]
Therefore, the symplectic form  $\Omega=\Omega_\cS^\epsilon$ on $\cS$ 
\bq
\Omega=
\begin{cases} 
\ \   p^*_-\omega^-  +  \pi_{-,\C}^*\ \omega_\C,                                                                            & \mathrm{on} \  \cE_-, \\  
\ \   \displaystyle\omega_{\cS_Z}^\epsilon=  \Pi_{Z}^*\omega_Z+\frac{\epsilon^2}{2} d(\alpha_\#),        & \mathrm{on} \  \cS_Z,  \\
\ \   p^*_+\omega^+  +  \pi_{+,\C}^*\ \omega_\C,                                                                         & \mathrm{on} \   \cE_+, \nonumber
\end{cases}
\eq
is $G$-invariant.  That is, $(\cS, \Omega)$ is a Hamiltonian $G$-manifold with  the associated moment map    $\mu:  \cS\to\g^*$  given by 
\bq\label{thm-MMP}
\mu=
\begin{cases} 
\ \   \mu^+ \circ p_+,                                                                            & \mathrm{on} \  \cE_-, \\  
\ \    \displaystyle\mu_Z\circ\Pi_{Z}+ \frac{\epsilon^2}{2} \mu_{\alpha_\#},        & \mathrm{on} \  \cS_Z,  \\
\ \    \mu^- \circ p_-,                                                                         & \mathrm{on} \   \cE_+. 
\end{cases}
\eq 

Assume that $0$ is a regular value of $\mu^\pm:  X^\pm\to \g^*$ and $G$ acts  freely on $(\mu^\pm)^{-1}(0)$.  Then  the quotient manifold,
the symplectic reductions of $(X^\pm, \omega^\pm, \mu^\pm)$, respectively, 
\[
X^\pm_0 = (\mu^\pm)^{-1}(0)/G
\]
has a unique symplectic structure   $\omega^\pm_0$  such that 
\[
\pi_\pm^* (\omega^\pm_0) = \omega^\pm|_{(\mu^\pm)^{-1}(0)}, 
\]
where $\pi_\pm:  (\mu^\pm)^{-1}(0)\to X^\pm_0$  is the quotient map.  These symplectic manifolds 
$(X^\pm_0, \omega^\pm_0)$ have a common symplectic submanifold 
\[
(Z_0 = \mu_Z^{-1}(0)/G, \omega_0^Z)
\]
with the opposite symplectic  normal bundles.  Here $\omega_0^Z$ is the unique symplectic structure on $Z_0$ satisfying
\[
(\pi_{Z_0})^* ( \omega_0^Z) =  \omega_Z |_{ \mu_Z^{-1}(0)},
\]
for the quotient map $\pi_{Z_0}: \mu_Z^{-1}(0)  \to \mu_Z^{-1}(0)/G$.    The 
$G$-equivariant isomorphism  (\ref{iso-1}) also descends to an   isomorphism of complex line bundles
\[
\Phi_0:N_{X_0^+} Z_0 \otimes_\C N_{X_0^-} Z_0 \longrightarrow  Z_0\times\C.
\]
So we can apply the symplectic sum operation to 
 $(X^+_0, \omega^+_0)$  and  $(X^-_0, \omega^-_0)$   
along  $(Z_0, \omega_0^Z)$
Note that  the symplectic reduction of $ ( \cS_Z, \omega^\epsilon_{\cS_Z}, \mu_{\cS_Z})$ is 
\[
\big((\mu_{\cS_Z})^{-1} (0)/G, ( \omega_{\cS_Z}^\epsilon)_0 \big)
\]
with  the unique symplectic form   $( \omega_{\cS_Z}^\epsilon)_0$ specified by 
 \[
( \pi_{\cS_Z, 0} )^* ( \omega_{\cS_Z}^\epsilon)_0 = \omega_{\cS_Z}^\epsilon |_{(\mu_{\cS_Z})^{-1} (0)}
\]
where $ \pi_{\cS_Z, 0}$ is the quotient map $(\mu_{\cS_Z})^{-1} (0) \to (\mu_{\cS_Z})^{-1} (0)/G$.  

From the symplectic sum construction of   $(X^+ , \omega^+)$  and  $(X^- , \omega^- )$  , we have
\[
\mu^{-1} (0) = (\mu_{\cE_-})^{-1}(0)   \  \bigcup_{gl_-}\ (\mu_{\cS_Z})^{-1} (0) \ \bigcup_{gl_+}\ (\mu_{\cE_+})^{-1}(0),
\]
where the gluing maps  $gl_\pm:   (\mu_{\cS_Z})^{-1} (0) \ \cap  \cS_{Z, \pm}  \to     (\mu_{\cE_\pm})^{-1}(0) $ are well-defined  open maps due to 
\[
 (\mu_{\cS_Z})^{-1} (0) \ \cap  \cS_{Z, \pm} = \big(\mu_{\cS_Z}|_{\cS_{Z, \pm}}\big)^{-1} (0)  
 \]
 and $\mu_{\cS_Z }|_{ \cS_{Z, \pm} }= \mu_{\cE_\pm} \circ gl_\pm $.   These gluing map descend to open maps
 \[
 \overline{gl}_\pm:  \big(\mu_{\cS_Z}|_{\cS_{Z, \pm}}\big)^{-1} (0) /G \to  (\mu_{\cE_\pm})^{-1}(0) /G
 \]
 of symplectic manifolds. 
  Therefore, 
 \[
 \cS\sslash G= \mu^{-1}(0) /G = (\mu_{\cE_-})^{-1}(0) /G   \  \bigcup_{\overline{gl}_-}\ (\mu_{\cS_Z})^{-1} (0)/G \ \bigcup_{\overline{gl}_+}\ (\mu_{\cE_+})^{-1}(0)/G,
 \]
 is a symplectic manifold with a unique  symplectic structure $\Omega_0$ given by
 \[
 \Omega_0=
\begin{cases} 
\ \   p^*_-\omega_0^- +  \pi_{-,\C}^*\ \omega_\C,                                                                            & \mathrm{on} \  (\mu_{\cE_-})^{-1}(0) /G , \\  
\ \   &  \\
\ \  \ \ \ \ \ \    (\omega_{\cS_Z}^\epsilon)_0,        & \mathrm{on} \  (\mu_{\cS_Z})^{-1} (0)/G,  \\
\ \   &  \\
\ \   p^*_+\omega_0^+  +  \pi_{+,\C}^*\ \omega_\C.                                                                    & \mathrm{on} \  (\mu_{\cE_+})^{-1}(0)/G,
\end{cases}
\]
Moreover,   there is a smooth map  
$\pi_0:  \cS\sslash G  \to D$ where $D$ is  a  sufficiently small neighbourhood of the origin in   $ \C$ such that 
 \begin{itemize}
\item  $\pi_0$ is surjective and  $\pi_0^{-1}(0) = X_0^+\cup_{Z_0}X_0^-$, 
\item $\pi_0$ is submersion away from $Z_0\subset \cS\sslash G$,
\item the restriction of $\Omega_0$ to $\pi_0^{-1}(\lambda)$ is nondegenerate for every $\lambda\in D\setminus \{0\}$,
\item $\Omega_0 |_{X_0^\pm} = \omega_0^\pm$ for $X_0^\pm \subset  \pi_0^{-1}(0)$.  
\end{itemize}
By the uniqueness of symplectic form for the symplectic reduction, for any $\lambda\in D-\{0\}$, the symplectic manifold
\[
(\pi_0^{-1}(\lambda), \Omega_0|_{\pi_0^{-1}(\lambda)}) 
\]
is the sympletic sum of $(X^+_0, \omega^+_0)$  and  $(X^-_0, \omega^-_0)$   
along  $(Z_0, \omega_0^Z)$.  This completes the proof  Theorem \ref{thm-1}.  

\br If  $0$ is a regular value of $\mu^\pm:  X^\pm\to \g^*$ and $\mu_Z:  Z\to \g^*$  but   $G$-actions is not free, 
then $(X_0^\pm, \omega_0^\pm)$ are symplectic orbifolds with a common symplectic suborbifold
  $(Z_0, \omega_0^Z)$ which is normalizable in the sense of \cite{VM}.  In this case, $(\cS\sslash G, \Omega_0)$ is a symplectic orbifold with a smooth map $\pi_0:  \cS\sslash G  \to D$ such that 
   for any $\lambda\in D-\{0\}$, the symplectic orbifold
\[
(\pi_0^{-1}(\lambda), \Omega_0|_{\pi_0^{-1}(\lambda)}) 
\]
is  a sympletic orbifold  sum of $(X^+_0, \omega^+_0)$  and  $(X^-_0, \omega^-_0)$   
along  $(Z_0, \omega_0^Z)$ as defined in \cite{VM}.  We leave the details of this construction to interested readers. 
  
\er

\br    For a  Hamiltonian $G$-manifold  with a local  Hamiltonian   $S^1$-action   which commutes with $G$-action,    the symplectic cut   produces  two    Hamiltonian $G$-manifolds (Remark 1.2  in \cite{L1}).  For completeness, we give an explicit description of this construction. Given a  Hamiltonian $G$-manifold  $(X, \omega, \mu)$  with a local  Hamiltonian   $S^1$-action   
which commutes with $G$-action. Let $V$ be  a $G$-invariant  open  subset   of $X$ with a $G$-commuting Hamiltonian $S^1$-action and  a $G$-invariant moment map 
\[
H:    V \longrightarrow  \R. 
\]
 Assume that  there is a small interval $I= (-\delta, \delta)$ of regular values and  $Y= H^{-1}(0)$ is a separating  compact hypersurface  with a free $S^1$-action.   By the  symplectic $S^1$-reduction,  there is a    circle bundle $\pi:  Y\to Z= Y/S^1$ and a  symplectic structure $\omega_Z$ on $Z$ 
uniquely defined by 
\[
\pi^*\omega_Z = \omega|_Y.
\]
As the $G$-action on $Y$ commutes with the $S^1$-action,   $\pi:  Y\to Z= Y/S^1$ is  a $G$-equivariant  circle bunlde and $(Z, \omega_Z)$ is a Hamiltonian $G$-manifold.  
Choosing  a  $G$-invariant connection $\alpha $ on $Y$, there exists   a  $G$-invariant 1-form  on $Y$,  also denoted by $\alpha$,     vanishing on
horizontal vector fields  defined by the connection and  
\bq\label{alp-Y}
\alpha  \left(\frac{d}{d\theta} e^{\i\theta} y \Big|_{\theta=0}\right)   =1,\ \ \  \forall\ y\in Y.
\eq
For simplicity,  identifying $V=H^{-1} (I) $ with $ I \times Y$, the symplectic form $\omega$  on $ I \times Y$  can be written as
\[
\omega=\pi^* \omega_Z + d(t \alpha)
\]
where $t$ is the coordinate on $I$. Moreover, let $\mu_\alpha: Y\to \g$ be the $\alpha$-moment map on $Y$, then restricted to $I\times Y$, we have 
\[
\mu = \mu_Z\circ \pi + t \mu_\alpha,
\]
where $\mu_Z: Z\to \g $ is  a moment map for the Hamiltonian $G$-action on $Z$.

 Let $\omega_\C$ denote the standard symplectic structure on $\C$.   Then  $(V   \times \C, \omega\oplus \omega_\C)$ be  a Hamiltonian $G\times S^1$-manifold, where 
 $G$ acts trivially on $\C$ and  the $S^1$-action on $\C$  is  given  by the scalar multiplication $e^{\pm \i \theta}$.  The  moment maps for 
 $S^1$-actions  are 
 \[
 H_\pm:  V \times \C \longrightarrow \R, \qquad H_\pm (x, z) = H(x) \mp \dfrac 12 |z|^2. 
 \]
   Let  $ (V_\pm = H_\pm ^{-1}(0) /S^1, \omega_0^\pm)$
 denote the $S^1$-symplectic reductions.  Then 
 \[
 (Z  =   \big(H_\pm^{-1}(0) \cap (V\times \{0\}) \big) /S^1, \omega_Z)
 \]
 is a $G$-invariant symplectic submanifold of $ V_\pm$ with  $G$-equivariant symplectic normal bundle  
 \[
 N_\pm  \cong Y\times_{\rho_\pm} \C.  
 \]
 Here  $\pi_\pm:  Y\times_{\rho_\pm} \C \to Z$ is the  complex line bundle over $Z$   associated  to the representations
 $ e^{\i \theta}\cdot z  = e^{\mp \i \theta}z$ for $z\in \C$.

 The $G$-invariant connection  $\alpha$ on $Y$ induces    $G$-invariant connection $\nabla^ \pm$  on $N_\pm$.  Choose $G$-invariant  Hermitian metrics  $g_\pm$ on $N_\pm$ which are  compatible with  $\nabla^ \pm$.  
With respect to these metrics, we can identify the $\delta$-neighbourhood of $Z$ in $V_\pm$ for a sufficiently small $\delta$,  denoted by
$N_\pm^\delta$, as
\[
N_\pm^\delta \cong \{[y, t, z]|  |z|^2 =\pm 2t, |t| \leq \delta/2\} \subset V_\pm
\]
with the induced symplectic form and the associated moment map
given by
\[
\omega^\delta_\pm = \pi^*_\pm \omega_Z\pm d (t \alpha), \qquad 
\mu_\pm^\delta = \mu_Z\circ\pi_\pm \pm t \mu_\alpha.
\]
 The symplectic cut is a union of two symplectic manifolds obtained by gluing $V_\pm$ and $X-Y$ through the following 
 symplectomorphisms
 \[
 \Psi_+: Y\times (0, \epsilon) \longrightarrow N_+^{2\epsilon}-Z, \qquad 
 \Psi_-: Y\times (-\epsilon, 0) \longrightarrow N_-^{2\epsilon}-Z
 \]
 defined by $\Psi_\pm (y, t) = [y, t,  \sqrt{\pm 2t}]$ for a sufficiently small $\epsilon$.  
  These result in two  Hamiltonian $G$-manifolds  
\bq\label{S-cut}
 X^+:= X^+_{0}\cup_{\Psi_+} V_+ \ \ \  \   \mathrm{and} \ \ \ \  
 X^-:= X^-_{0}\cup_{\Psi_+} V_-, 
 \eq
 where $X_0^+$ and $X_0^-$ are  two components of $X-Y$ with end modelled on
 $Y\times (0, \epsilon)$ and   $Y\times (-\epsilon, 0)$ respectively. 
 The normal bundles  $N_{\pm}  $  of $Z$ in $X^\pm$ have opposite equivariant Euler classes.  
 
\er

\section{Equivariant first Chern classes}\label{equi-1-Chern}
 
In this section,  we compare the equivariant first Chern class of  the Hamiltonian sum  $(\cS_\lambda, \Omega_\lambda, \mu_\lambda)$  with   the equivariant first Chern classes of $(X^+,\omega^+,  \mu^+)$ and     $(X^-,\omega^-,  \mu^-)$.    Let $J_\lambda$ be  the $G$-invariant  $ \Omega_\lambda$-compatible  almost complex structure on $\cS_\lambda$ induced from the $G$-invariant  $ \Omega$-compatible  almost complex structure on $\cS$, and  $J^\pm$ be  the  induced  $G$-invariant  $\omega^\pm$-compatible almost complex structure on $X^\pm$. 
Denote by
\[
c_1^G(TX^\pm) \in H^2_G (X^\pm_G, \Z), \qquad c_1^G(T\cS_\lambda) \in H^2_G ((\cS_\lambda)_G, \Z)
\]
the first equivariant Chern classes  of $TX^\pm$ and $T\cS_\lambda$.  Here   $X^\pm_G = EG\times_G X^\pm$ and 
$(\cS_\lambda)_G = EG \times_G \cS_\lambda$ are the homotopy quotients of $X^\pm$ and $\cS_\lambda$ respectively. 

 We claim that  $\cS_\lambda$  can be realised as $X^+\#_{\varphi_\lambda} X^-$ obtained from gluing the complements of tubular neighbourhoods of  $Z$  in $X^+$  and  $X^-$  along their  boundaries   by an orientation-reversing $G$-equivariant
diffeomorphism  $\varphi_\lambda$ (defined by $\lambda$ and  the  $G$-equivariant Hermitian  line bundle  isomorphism (\ref{iso-1})).   To see this, recall that
\bq\label{cS_lam} 
\cS_\lambda = \big(X^- - \Psi_{-;\epsilon}(\overline{N_{-}Z(1)}) \big) \ \bigcup_{gl_-}\ \cS_{Z, \lambda} \ \bigcup_{gl_+}\ \big(X^+ - \Psi_{+;\epsilon}(\overline{N_{+}Z(1)} ) \big),
\eq
where the middle  part $\cS_{Z, \lambda}  =  \{(z,v,w)\in N_{+}Z\oplus N_{-}Z\ \big|\ |v|,|w|< 2, \epsilon \cP(v,w) =\lambda \}$ can be identified  as  the  open annular bundle  of  $N_+Z$ with radius from $\epsilon^{-1}  |\lambda|/2$ to $2$
\[
\big(  N_{+}Z(2) - \overline{N_{+}Z(\epsilon^{-1}  |\lambda|/2)} \big), 
\]
as we can  solve   $w$ in terms of $v$ from the equation $ \epsilon \cP(v,w) =\lambda $.    From the construction of the symplectic sum, we have
\[
\epsilon^{-1} |\lambda|  < \epsilon^{-1} \delta < \frac 12   \Longrightarrow  \epsilon^{-1}  |\lambda|/2 < \frac 14. 
\]
Therefore  $\cS_{Z, \lambda} $ consists of  three parts:
 \begin{itemize}
\item  the domain of  the  gluing maps $gl_+$ 
\begin{eqnarray} 
\cS_{Z, \lambda, +}  & = &   \{(z,v,w)\in N_{+}Z\oplus N_{-}Z\ \big| 1< |v| < 2, \epsilon \cP(v,w) =\lambda \} \nonumber\\ 
  & \cong&   \big(  N_{+}Z(2) - \overline{N_{+}Z(1)}\big),   \nonumber
\end{eqnarray} 
under this identification, $gl_+ = \Psi_{+;\epsilon}.$ 
\item the neck part  $\big(   \overline{N_{+}Z(1) } -  N_{+}Z(\epsilon^{-1}  |\lambda|) \big)$,   
which is  a closed annular bundle  of  $N_+Z$ with radius from  $\epsilon^{-1} |\lambda|$ to $1$. Note that  the neck part is identified with  $ \Psi_{+;\epsilon}(\ \overline{N_{+}Z(1) } -  N_{+}Z(\epsilon^{-1}  |\lambda|) )$, 
\item the domain of  the gluing map  $gl_-$  
\begin{eqnarray} 
 \cS_{Z, \lambda, -}  & = &   \{(z,v,w)\in N_{+}Z\oplus N_{-}Z\ \big| 1< |w| < 2, \epsilon \cP(v,w) =\lambda \}  \nonumber\\ 
  & \cong& \big(  N_{+}Z(\epsilon^{-1} |\lambda|) - \overline{N_{+}Z(\epsilon^{-1} |\lambda|/2)}\big).   \nonumber
\end{eqnarray}  
\end{itemize}
Thus $\cS_\lambda$ as in (\ref{cS_lam}) can  be written as
\[
\cS_\lambda = \big(X^- - \Psi_{-;\epsilon}(\overline{N_{-}Z(1)}) \big) \ \bigcup_{\varphi_\lambda}\  
 \big(X^+ - \Psi_{+;\epsilon}(\overline{N_{+}Z(\epsilon^{-1} |\lambda|)} ) \big),
\]
where $\varphi_\lambda: \partial \big(X^+ - \Psi_{+;\epsilon}(\overline{N_{+}Z(\epsilon^{-1} |\lambda|)} ) \big) \to  \partial \big(X^- - \Psi_{-;\epsilon}(\overline{N_{-}Z(1)}) \big)$ 
is the orientation-reversing $G$-equivariant diffeomorphism  defined by the   composition of the following three maps:
\begin{itemize}
\item  the restriction of gluing map
$ (gl_+)^{-1} = (\Psi_{+;\epsilon}.)^{-1}:  \partial \big(X^+ - \Psi_{+;\epsilon}(\overline{N_{+}Z(\epsilon^{-1} |\lambda|)} ) \big) \to SN_+Z (\epsilon^{-1} |\lambda|)$ (the circle bundle of $N_+Z$ of radius $\epsilon^{-1} |\lambda|$), 
\item an  orientation-reversing $G$-equivariant diffeomorphism  from $ SN_+Z (\epsilon^{-1} |\lambda|)$ to 
$SN_-Z(1)$ (the unit circle bundle of $N_-Z$)  by solving  $w$ in terms of $v\in SN_+Z (\epsilon^{-1} |\lambda|)$  from  the equation
\[
\epsilon \cP (v, w) = \lambda,
\]
\item  the restriction of gluing map $gl_-:   SN_-Z(1)  \to   \partial \big(X^- - \Psi_{-;\epsilon}(\overline{N_{-}Z(1)}) \big)$. 
\end{itemize}

 For the comparison,  we evaluate  
$c_1^G(TX^\pm)$ on  certain  degree 2  equivariant homology  classes  of the form
\[
B_\pm = (u^\pm_G)_* ([\Sigma_\pm])  \in H_2^G(X^\pm_G, \Z) 
\]
 coming from a principal $G$-bundle $P^\pm$ over   oriented closed  surfaces $\Sigma^\pm$ together with  smooth sections $u_\pm$ of $P^\pm\times_G X^\pm$ as in the following diagram
\bq\label{sym-vortex}
\xymatrix{
P^\pm\times_G X^\pm \ar[d]  \ar[r] & X^\pm_G \ar[d]   \\
\Sigma_\pm  \ar@/^2.0pc/[u]^{u^\pm}  \ar[r]_{c_P} \ar@{-->}[ru]^{u^\pm_G} &  BG,  } 
\eq
where $c_P:  \Sigma_\pm  \to BG$ is the classifying map of $P^\pm$.

 We assume that sections $u^\pm$  in (\ref{sym-vortex})  satisfy the following conditions
\bq\label{u-pm1}
(u^\pm)^{-1} (P^\pm \times_G Z) = \{x_1^\pm, \cdots, x_k^\pm\} \subset \Sigma_\pm, 
\eq
\bq\label{u-pm2}
u^+ (x_i^+) = u^-(x_i^-) \in P\times_G Z, \qquad \text{ord}_{x_i^+, Z}^{P^+} (u^+)  = \text{ord}_{x_i^-, Z}^{P^-} (u^-), 
\eq
for any $i =1, \cdots, k$. Here $\text{ord}_{x_i^\pm, Z}^{P^\pm}( u^\pm)$ is the order of contact of $u^\pm$ with $Z$ at $x_i^\pm$ defined as follows.  Take a local trivialization of $P^+$ in a small  neighborhood  $B_{x_i^+}$ of $x_i^+$ over which the section $u^+$ can be expressed as a smooth map from 
$B_{x_i^+}$ to $X^+$.    Take a  homotopy of 
$u^+$ on a small neighborhood of $ u^+(x_i^+)$ without changing the intersection  point $ u^+(x_i^+)$ so that 
$du^+$ maps a small circle  in  $T_{x_i^+} \Sigma^+$  to   a  circle   in $ N_+Z|_{u^+(x_i^+)}$.  
The degree of this map is  the order of contact  $\text{ord}_{x_i^+, Z}^{P^+} (u^+) \in \Z$.  
This degree is independent of the choice of local trivialization of $P^+$ near $x_i^+$ and the choice of the  homotopy of $u^+$ on a small neighborhood of $x_i^+$.    
The order of contact of $u^-$ with $Z$ at $x_i^-$ is defined in the same way.  

We can  construct a smooth section $u^+\#_{\varphi_\lambda}  u^-$ of $(P^+\#_\lambda P^-)\times_G \cS_\lambda$   for  a principal $G$-bundle  
$P^+\#_\lambda P^-$ over  a $\Sigma_+\#_\lambda\Sigma_-$ by the following steps
\begin{itemize}
\item removing   small discs  $D_{x_i^+}$ and $D_{x_i^-}$ around   $x_i^\pm$ in $\Sigma^\pm$  for all $i$ to get surfaces  $\widehat{\Sigma}^\pm$ with $k$ boundary circles. We can assume (by homotopying if necessary) 
\[
u^\pm|_{\partial D_{x_i^\pm}}:  \partial D_i \longrightarrow \Psi_{\pm, \epsilon} \big( \partial N_\pm Z (1)_{u^\pm (x_i^\pm)}\big), 
\]
\item gluing $\widehat{\Sigma}^+$ and $\widehat{\Sigma}^-$ along  each of  boundary circles by   orientation-reversing diffeomorphisms
$\{\varphi_i \}$    to form a smooth oriented  surface $\Sigma^+\#_\lambda \Sigma^-$   with neck structure
   such that  
\[
\Sigma^+\#_\lambda \Sigma^- =  \widehat{\Sigma}^-  \cup_{\{\varphi_i\}}   \big(  \bigsqcup_{i=1}^k  [\epsilon^{-1}  |\lambda|, 1]_i\times S^1 \big) \cup \widehat{\Sigma}^+.
\]

\item gluing the principal bundles $P^+$ and $P^-$ along each  boundary circle to form a principal bundle  over $\Sigma^+\#_\lambda \Sigma^-$, 
by using the trivializations $P^\pm|_{D_{x_i^\pm}} \cong D_{x_i^\pm} \times G$ and gluing parameters 
$ \rho_i: P^+ |_{x_i^+} \longrightarrow  P^-|_{x_i^-} $ ($G$-equivariant diffeomorphisms) for all $i$.  
As the isomorphism class of  the resulting principal  bundle doesn't depend on the gluing parameter, we will suppress  the gluing parameters and simply denote this bundle by
 $P^+  \#_\lambda  P^-$.  
 \item  homotopying $u^+$   on a small neighbourhood of $\partial D_{x_i^+}$   within a  small ball  $B_{u^+(x^+_i)}$  around 
 $u^+(x^+_i)$  in $X^+ $     for each  $i$ so that 
 $u^-\circ \varphi_i = \varphi_\lambda\circ u^+$ on  $\partial D_{x_i^+}$   for all $i$. To be precise,   choose  a homotopy
 \[
\tilde u^+  =  \{u_t^+\}_{t\in [[\epsilon^{-1}  |\lambda|, 1]} :    \bigsqcup_{i=1}^k  [\epsilon^{-1}  |\lambda|, 1]_i\times S^1 \to   \bigsqcup_{i=1}^k gl_+   \big(   \overline{N_{+}Z(1) } -  N_{+}Z(\epsilon^{-1}  |\lambda|) \big)_{u^+ (x_i^+)}
 \]
 such that  for each $i$,   $u_t^+|_{t=1}  = u^+|_{\partial D_{x_i^+}}$   and $u^-|_{\partial D_{x_i^-}}\circ \varphi_i = \varphi_\lambda\circ u_t^+|_{t=\epsilon^{-1}  |\lambda|}$.     This can be achieved because $P^\pm\times_G X^\pm$ are trivialized over $D_{x_i^\pm}$ , so $u^\pm$ can be written as smooth maps  $\partial D_{x_i^\pm} \to X^\pm$ and the degrees of maps
$\varphi_\lambda\circ u^+\circ \varphi_i^{-1}$ and $u^-$ are the same. 
\item We remark that $\cS_\lambda$ contains a neck part  $gl_+   \big(   \overline{N_{+}Z(1) } -  N_{+}Z(\epsilon^{-1}  |\lambda|) \big)$.    Define 
\bq\label{u-sum}
u^+\#_{\varphi_\lambda}  u^- =
\begin{cases} 
\ \    u^-  & \mathrm{on} \  \widehat{\Sigma}^-, \\ 
\ \   \tilde u^+  & \mathrm{on} \     \bigsqcup_{i=1}^k  [\epsilon^{-1}  |\lambda|, 1]_i\times S^1  \\   
\ \    u^+ & \mathrm{on} \  \widehat{\Sigma}^+. 
\end{cases}
\eq
Then $u^+\#_{\varphi_\lambda}  u^- $ is a section of $(P^+  \#_\lambda  P^-)\times_G \cS_\lambda$.
\end{itemize}
From  the construction of $u^+\#_{\varphi_\lambda}  u^-$ and the homotopy invariance of homology theory,   $u^+$ and $u^-$ completely define the homology class of 
\[
(u^+\#_{\varphi_\lambda}  u^- )_* ([\Sigma^+\#_\lambda \Sigma^-]) \in H_2(P\times_G \cS_\lambda,  \Z),
\]
thus an equivariant homology classes $\big((u^+\#_{\varphi_\lambda}   u^-)_G\big)_* \in H_2^G(\cS_\lambda, \Z)$.  

The Hamiltonian sum construction gives rise to a Hamiltonian $G$-manifold of the form $X^+\#_{\varphi_\lambda} X^-$.   We can compare equivariant first Chern class of  the Hamiltonian sum with the equivariant first Chern classes of 
$X^+$ and $X^-$.  

\bp   Suppose that the sections $u^\pm$  in (\ref{sym-vortex})  satisfy the   conditions (\ref{u-pm1}) and (\ref{u-pm2}). Let  $B_\pm$ be the equivariant homology classes $ (u^\pm_G)_* ([\Sigma_\pm]) \in H_2^G(X^\pm, \Z)$   and $B_+\#_{\varphi_\lambda} B_-$ be the equivariant homology class in $H_2^G(\cS_\lambda, \Z)$   defined by    $u^+ \#_{\varphi_\lambda} u^-$, then 
\[
\la c_1^G(T(\cS_\lambda )),B_+\#_{\varphi_\lambda} B_ - \ra     
= \la c_1^G(TX^+), B_+ \ra  + \la c_1^G(TX^-), B_-\ra  -  2 \sum_{i=1}^k \text{ord}_{x_i^+, Z}^{P^+} (u^+) .
\]
and $\ \sum_{i=1}^k \text{ord}_{x_i^+, Z}^{P^+} (u^+)  =  \sum_{i=1}^k \text{ord}_{x_i^-, Z}^{P^-} (u^-) $ is the intersection number  
\[
(u^+)_* ([\Sigma_+]) \cdot [P^+\times_G Z] = (u^-)_* ([\Sigma_-]) \cdot [P^-\times_G Z].
\]
\ep
 
 \begin{proof}  From the functoriality of the first Chern class, we have 
\[
\la c_1^G(TX^\pm),  (u^\pm_G)_* ([\Sigma_\pm]) \ra = \la c_1 (P\times_G TX^\pm),  (u^\pm)_*([\Sigma_\pm])\ra,
\]
and 
\[
\la c_1^G(T\cS_\lambda), ((u^+\#_{\varphi_\lambda}  u^-)_G)_* ([\Sigma^+\#_\lambda \Sigma^-]) \ra = \la c_1 (P\times_G T\cS_\lambda),  (u^+\#_{\varphi_\lambda}  u^-)_*([\Sigma^+\#_\lambda \Sigma^-] )\ra,
\]
In $H_2^G (\cS, \Z)$, the  equivariant homology class  of $B_+\#_{\varphi_\lambda} B_-$ is the sum of   equivariant homology class  of  
$ B_+$ and $ B_-$. Note that the normal bundle of  $\cS_\lambda$ in $\cS$ is trivial, we have 
\ba\label{c1-sum}
&& \la c_1 (P\times_G T\cS_\lambda),  (u^+\#_{\varphi_\lambda}  u^-)_*([\Sigma^+\#_\lambda \Sigma^-] )\ra  \nonumber\\
& = &  \la c_1 (P\times_G T\cS ),  (u^+)_* [\Sigma^+] \ra  +    \la c_1 (P\times_G T\cS ),  (u^-)_* [\Sigma^-] \ra . 
\na
Note that  $u^\pm$ are sections of $P\times_G X^\pm \subset P\times_G  \cS_0$. Along $X^+\subset \cS_0$,  we have  the decomposition of $G$-equivariant vector bundle, 
\[
T\cS|_{X^+}  \cong  TX^+ \oplus N_\cS X^+.
\]
As  a  $G$-equivariant vector bundle,  
 $N_\cS X^+$ is trivial on $X^+- \Psi_{+;\epsilon}(\overline{N_{+}Z(1)})$  and isomorphic to  $(\Psi_{+;\epsilon}^{-1})^* \circ \pi_{+, Z}^* N_-Z$  on the   $\Psi_{+;\epsilon}( N_{+}Z(2) )$,  where 
 $\pi_{+, Z}:  N_{+}Z(2) \to Z$ is the restriction of the bundle projection.  As $N_-Z$ is complex  conjugate to $N_+Z$,  we 
 get
 \[
 c_1(  P^+\times_G  N_\cS X^+) = - PD ([P^+\times_G Z]) \in H^2(P^+\times_G X^+, \Z).
 \]
 Thus 
 \[
  \la PD ([P^+\times_G Z]),   (u^+)_* [\Sigma^+] \ra =   (u^+)_* [\Sigma^+]\cdot  [P^+\times_G Z]   =  \sum_{i=1}^k \text{ord}_{x_i^+, Z}^{P^+} (u^+)
 \]
 Putting these together, we have 
 \[ 
  \la c_1 (P\times_G T\cS ),  (u^+)_* [\Sigma^+] \ra  =  \la c_1 (P\times_G TX^+ ),  (u^+)_* [\Sigma^+] \ra 
    -    \sum_{i=1}^k \text{ord}_{x_i^+, Z}^{P^+} (u^+).
      \] 
Similarly,  we have 
  \[
  \la c_1 (P\times_G T\cS ),  (u^-)_* [\Sigma^-] \ra  =  \la c_1 (P\times_G TX^- ),  (u^-)_* [\Sigma^-] \ra 
 -    \sum_{i=1}^k \text{ord}_{x_i^-, Z}^{P^-} (u^-). 
      \]
 Together with  (\ref{c1-sum}), this completes the proof of the Proposition.

 \end{proof}

\bigskip
\noindent\small School of Mathematics\\
\noindent\small Sichuan University, Chengdu 610065, China\\ 
 \noindent\small E-mail address: chenbohui@scu.edu.cn

\bigskip
\noindent\small Department of Mathematics\\
\noindent\small Jinan University, Guangzhou 510632, China\\
\noindent\small E-mail address:  hailongher@jnu.edu.cn

\bigskip
\noindent\small Mathematical Sciences Institute\\
\noindent\small The National Australian University, Canberra ACT 2601, Australia\\
\noindent\small E-mail address:  bai-ling.wang@anu.edu.au

\end{document}